\documentclass[11pt]{article}
\usepackage{amssymb}
\newtheorem{theorem}{THEOREM}[section]
\newtheorem{proposition}{PROPOSITION}[section]
\newtheorem{lemma}{LEMMA}[section]
\newtheorem{corollary}{COROLLARY}[section]
\newtheorem{remark}{REMARK}[section]

\begin{document}
\title{The Best Constant, the Existence of
Extremal Functions and Related Results for an Improved Hardy-Sobolev
Inequality}
\author{N. B. Zographopoulos \\ \\ Department of Science, Division
of Mathematics,
\\
Technical University of Crete, \\
73100 Chania, Greece\\
E-mail address: \ \  nzogr@science.tuc.gr }
\date{}
\maketitle
\pagestyle{myheadings} \thispagestyle{plain} \markboth{N. B.
Zographopoulos}{}
\maketitle
\begin{abstract}
We present the best constant and the
existence of extremal functions for an Improved Hardy-Sobolev
inequality. We prove that, under a proper transformation, this
inequality is equivalent to the Sobolev inequality in
$\mathbb{R}^N$. We also discuss the connection of the related
functional spaces and as a result we obtain some Caffarelli - Kohn
- Nirenberg inequalities. Our starting point is the existence of a
minimizer for the Bliss' inequality and the indirect dependence of
the Hardy inequality at the origin.
\end{abstract}
 \vspace{0.1cm}
%
%
%
%
\section{Introduction}
Assume the following inequality:
\begin{equation}\label{m}
\int_{0}^{R} r\, |v'|^2\, dr \geq c\, \left ( \int_{0}^{R}
r^{-1}\, \left ( - \log \left ( \frac{r}{R} \right ) \right
)^{-\frac{2(N-1)}{N-2}} |v|^\frac{2N}{N-2}\, dr \right
)^{\frac{N-2}{N}},
\end{equation}
which holds for any function $v \in C_{0}^{\infty} (0,R)$. This
inequality may be obtained from a more general inequality
\cite[Theorem 4]{maz85} (see also \cite[Lemma 2.2]{ft02})).
However, as prof. V. Maz'ya pointed to us this inequality is also
obtained from Bliss' inequality \cite{bliss} ( For the derivation
of this inequality and some related discussion we refer to Section
3). In this work, we prove that under a proper transformation
inequality (\ref{m}) is equivalent to the Sobolev inequality in
$\mathbb{R}^N$ and consequently we obtain the best constants and
the minimizers for (\ref{m}). The best constant in the Sobolev
inequality in $\mathbb{R}^N$:
\begin{equation}\label{sobolev}
\int_{\mathbb{R}^N} |\nabla u|^2 \, dx \geq S\, \left (
\int_{\mathbb{R}^N} |u|^{\frac{2N}{N-2}}\, dx \right
)^{\frac{N-2}{N}},
\end{equation}
as it is well known, see \cite{Au,LiThi,Tal}, is
\[
S(N) = \frac{N(N-2)}{4}\, |\mathbb{S}_N|^{2/N} = 2^{2/N}\,
\pi^{1+1/N}\, \Gamma \left ( \frac{N+1}{2} \right )^{-2/N},
\]
where $\mathbb{S}_N$ is the area of the N-dimensional unit sphere
and the extremal functions are
\[
\psi_{\mu,\nu}(|x|) = (\mu^2 + \nu^2 |x|)^2)^{-(N-2)/2},\;\;\; \mu
\ne 0,\, \nu \ne 0.
\]
For a quantitative version of the sharp Sobolev inequality we
refer to \cite{cfmt}.
\begin{lemma} \label{lemmam}
Inequality (\ref{m}) under the transformation
\begin{equation}\label{trans}
 u(r) = w(t),\;\;\; t = \left ( - \log \left ( \frac{r}{R} \right )
 \right
)^{-\frac{1}{N-2}}
\end{equation}
is equivalent to (\ref{sobolev}). The best constant is
\begin{equation} \label{bestm}
C_M= (N-2)^{-\frac{2(N-1)}{N}} (N\, \omega_N)^{-\frac{2}{N}}\,
S(N) = \frac{1}{4}\, \left ( \frac{N}{N-2} \right
)^{\frac{N-2}{N}}\, \left ( \frac{|\mathbb{S}_N|}{\omega_N} \right
)^{\frac{2}{N}},
\end{equation}
where $\omega_N$ denotes the Lebesgue measure of the unit ball in
$\mathbb{R}^N$ and the minimizers are
\begin{equation} \label{minimm}
\phi_{\mu,\nu}(r) = \psi_{\mu,\nu} (t) = \left ( \mu^2 + \nu^2
\left ( - \log \left ( \frac{r}{R} \right ) \right
)^{-\frac{2}{N-2}} \right )^{-\frac{N-2}{2}},\;\;\; \mu \ne 0,\,
\nu \ne 0.
\end{equation}
\end{lemma}
It is clear that $\phi$ may be continuously defined as
$\phi_{\mu,\nu}(0) = \mu^{-(N-2)}$ and $\phi_{\mu,\nu}(R)=0$.
\vspace{0.2cm}

As an application of inequality (\ref{m}) the authors in
\cite{ft02}, proved the following Improved Hardy-Sobolev
(IHS) inequality:
\begin{eqnarray}\label{HS1}
\int_{B_R} |\nabla u(|x|)|^2 \, dx &\geq& \left ( \frac{N-2}{2}
\right )^2 \int_{B_R} \frac{u^2(|x|)}{|x|^2} \, dx \nonumber \\
&& + C_{HS} \left ( \int_{B_R} |u(|x|)|^{\frac{2N}{N-2}}
\left(-\log \left ( \frac{|x|}{R} \right)
\right)^{-\frac{2(N-1)}{N-2}}\, dx \right )^{\frac{N-2}{N}},
\end{eqnarray}
in the radial case, i.e. where $B_R$ is the open ball in
$\mathbb{R}^N$, $N \geq 3$, of radius $R$ centered at the origin
and $u \in C_{0}^{\infty} (B_R \backslash \{0\})$ is a radially
symmetric function. The same result was proved in \cite{mus}, with
the use of a Caffarelli-Kohn-Nirenberg inequality. Actually, in
\cite{ft02} the following general (not in necessarily radial case)
IHS inequality was proved: Let $\Omega$ be a bounded domain in
$\mathbb{R}^N$, $N \geq 3$, containing the origin, $D_0 = \sup_{x
\in \Omega} |x|$ and $D>D_0$, then the following inequality
\begin{eqnarray} \label{HS}
\int_{\Omega} |\nabla u|^2 \, dx &\geq& \left ( \frac{N-2}{2}
\right )^2 \int_{\Omega} \frac{u^2}{|x|^2} \, dx \nonumber \\ &&+
C_{HS} (\Omega)\, \left ( \int_{\Omega} |u|^{\frac{2N}{N-2}}
\left(-\log \left ( \frac{|x|}{D} \right)
\right)^{-\frac{2(N-1)}{N-2}}\, dx \right )^{\frac{N-2}{N}}
\end{eqnarray}
holds for any $u \in C_0^{\infty} (\Omega \backslash \{0\})$. From
the discussion in \cite{ft02, mus}, it is clear that the nature of
(\ref{HS}) depends on the distance of $D$ from $D_0$, for instance
in the case where $D = D_0$ the author in \cite{mus} proved that
the inequality cannot hold if we consider nonradial functions.

Both papers follow the approach that is based on the following change of
variables (This approach was introduced in \cite{bv} and followed in various ways
by many authors); For any $u \in H_0^1 (\Omega)$ we set
\begin{equation} \label{changevarv}
u = |x|^{\frac{N-2}{2}}\, v
\end{equation}
and in this case we have
\begin{equation}\label{W=H}
\int_{\Omega} |\nabla u|^2 \, dx - \left ( \frac{N-2}{2} \right
)^2 \int_{\Omega} \frac{u^2}{|x|^2} \, dx = \int_{\Omega}
|x|^{-(N-2)}\, |\nabla v|^2 \, dx.
\end{equation}
Then, inequality (\ref{HS}) is equivalent to
\begin{equation}\label{HSv}
\int_{\Omega} |x|^{-(N-2)}\, |\nabla v|^2 \, dx \geq C_{HS}
(\Omega)\, \left ( \int_{\Omega} |x|^{-N} |v|^{\frac{2N}{N-2}}
\left(-\log \left ( \frac{|x|}{D} \right)
\right)^{-\frac{2(N-1)}{N-2}}\, dx \right )^{\frac{N-2}{N}}
\end{equation}
and which in turn is equivalent, in the radial case, to (\ref{m}).
Therefore, is natural to consider the space $W_{0}^{1,2}
(|x|^{-(N-2)},\Omega)$, see \cite{ft02}, which is defined as the complement of the
$C_0^{\infty}(\Omega)$ functions under the norm
\[
||v||^2_{W_{0}^{1,2} (|x|^{-(N-2)},\Omega)} = \int_{\Omega}
|x|^{-(N-2)}\, |\nabla v|^2 \, dx.
\]
The space $W_{0}^{1,2} (|x|^{-(N-2)},\Omega)$ has the property
(for a generalization see Lemma \ref{lemmaWH}) that if $u \in
H_0^1 (\Omega)$ then $|x|^{(N-2)/2} u \in W_{0}^{1,2}
(|x|^{-(N-2)},\Omega)$. (Some other properties of this space may
be found in Section 2). The advantage of this space is the
following; assume an inequality, e.g. an improved Hardy inequality
(see \cite[Section 3]{ft02}), which admits no $H_0^1$- minimizer
then, under the change of variables (\ref{changevarv}), the
corresponding inequality admits $W_{0}^{1,2}
(|x|^{-(N-2)},\Omega)$ minimizer. This happens because if
$v \in W_{0}^{1,2} (|x|^{-(N-2)},\Omega)$
then it is not necessary that $|x|^{-(N-2)/2}\, v$ belongs in
$H_0^1 (\Omega)$.

From (\ref{W=H}) it is also natural to define the space $H$ as
the completion of the set
\[
\left \{ |x|^{-\frac{N-2}{2}}\, \phi (x); \phi \in
C_{0}^{\infty} (\Omega) \right \}
\]
under the norm
\begin{equation} \label{norm}
||u||^{2}_{H(\Omega)} = \int_{\Omega} |\nabla u|^2\, dx -
\left ( \frac{N-2}{2} \right )^2 \int_{\Omega}
\frac{u^2}{|x|^2}\, dx - L^2(u)
\end{equation}
where $u_r$ is the radial part of $u$, i.e. we extend $u$ as zero outside $\Omega$, and
for some $R > \sup_{x \in \Omega} |x|$, we take the projection of $u$ on the space of
radially symmetric functions, i.e.,
\[
u_r (|x|) = \frac{1}{|\partial B_R|} \int_{\partial B_R} u\, ds
\]
and by $L(u)$ we denote the quantity
\begin{equation} \label{L}
L(u) : = \left ( \frac{N(N-2)}{2}\, \omega_N \right)^{1/2}\, \lim_{|x| \to 0} |x|^{\frac{N-2}{2}} u_r (|x|).
\end{equation}
For the definition of this space and some related properties we
refer to \cite{vz, vz00}. We note that $H_0^1 (\Omega)$ is a subspace
of $H(\Omega)$. The fact that the space $H$ is not convenient
to be defined as the completion of the
$C_{0}^{\infty} (\Omega)$ functions under the norm
\begin{equation} \label{oldnorm}
||\phi||^{2}_{H(\Omega)} = \int_{\Omega} |\nabla \phi|^2\, dx -
\left ( \frac{N-2}{2} \right )^2 \int_{\Omega}
\frac{\phi^2}{|x|^2}\, dx.
\end{equation}
is explained in \cite{vz} and this due to the presence of a ``boundary'' term; if
we define $H$ with norm given by (\ref{oldnorm}) then functions that behave at
the origin like $|x|^{-(N-2)/2}$ fail to be in $H$.

The connection between the spaces $H(\Omega)$ and $W_{0}^{1,2}
(|x|^{-(N-2)},\Omega)$ is given in the following lemma;
\begin{lemma} \label{lemmaWH}
Assume that $\Omega$ is a bounded domain of $\mathbb{R}^N$, $N
\geq 3$, containing the origin. Then, $u \in H(\Omega)$ if and only if $|x|^{(N-2)/2}\, u \in W_{0}^{1,2}
(|x|^{-(N-2)},\Omega)$. In this case the connection of the norms is given by
\[
||u||^2_{H(\Omega)} = ||v||^2_{W_{0}^{1,2} (|x|^{-(N-2)},\Omega)},
\]
or
\begin{equation}\label{normWH}
\int_{\Omega} |\nabla u|^2 \, dx - \left ( \frac{N-2}{2} \right
)^2 \int_{\Omega} \frac{u^2}{|x|^2} \, dx - L^2(u)= \int_{\Omega}
|x|^{-(N-2)}\, |\nabla v|^2 \, dx.
\end{equation}
\end{lemma}

In addition we can relate these spaces, in the radial case, with
the space $D^{1,2}(\mathbb{R}^N)$, which is defined as the closure
of $C_{0}^{\infty}(\mathbb{R}^N)$ functions under the norm
\[
||\phi||^2_{D^{1,2}(\mathbb{R}^N)} = \int_{\mathbb{R}^N} |\nabla
\phi|^2\, dx.
\]
For more details we refer to the classical book \cite{ada75}. If
we denote by $H_r (\Omega)$, $W_{0,r}^{1,2} (|x|^{-(N-2)},\Omega)$
and $D_r^{1,2}(\mathbb{R}^N)$ the subspaces of $H (\Omega)$,
$W_{0}^{1,2} (|x|^{-(N-2)},\Omega)$ and $D^{1,2}(\mathbb{R}^N)$,
respectively, which consist of radial functions, we have that
\begin{lemma} \label{W=D}
Let $v \in W_{0,r}^{1,2} (|x|^{-(N-2)},B_R)$ and set
\begin{equation} \label{transxw}
v(|x|) = w(t),\;\;\;\;\;\; t = \left ( - \log \left (
\frac{|x|}{R} \right ) \right )^{-\frac{1}{N-2}}
\end{equation}
as in (\ref{trans}). Then, $v \in W_{0,r}^{1,2}
(|x|^{-(N-2)},B_R)$ if and only if $w \in D_r^{1,2}(\mathbb{R}^N)$
and
\begin{equation}\label{eqWD}
||v||_{W_{0,r}^{1,2} (|x|^{-(N-2)},B_R)} = (N-2)^{-1}\,
||w||_{D_r^{1,2}(\mathbb{R}^N)}.
\end{equation}
\end{lemma}
Observe that (\ref{eqWD}) is independent of the radius $R$ and in
the case where $N=3$ the norm in $W_{0,r}^{1,2}
(|x|^{-(N-2)},B_R)$ coincides with the norm in
$D_r^{1,2}(\mathbb{R}^N)$. Moreover, (\ref{normWH}) and
(\ref{eqWD}) imply that
\begin{corollary}
Let $u \in H_r (B_R)$ and set
\begin{equation} \label{transx}
w(t) = |x|^{\frac{N-2}{2}} u(|x|), \;\;\;\;\;\; t = \left ( - \log
\left ( \frac{|x|}{R} \right ) \right )^{-\frac{1}{N-2}}.
\end{equation}
Then, if $u \in H_r (B_R)$ then $w \in
D_r^{1,2}(\mathbb{R}^N)$ and
\begin{equation}\label{eqHD}
||u||^2_{H_r (B_R)} = (N-2)^{-1}\,
||w||^2_{D_r^{1,2}(\mathbb{R}^N)}.
\end{equation}
\end{corollary} \vspace{0.2cm}

For a related to the IHS inequality (\ref{HS1}), as a
consequence of Lemma \ref{lemmam}, we have
\begin{theorem} \label{thmHSw}
The infimum of the ratio
\begin{equation}\label{HSw}
\frac{\int_{B_R} |x|^{-(N-2)}\, |\nabla v(|x|)|^2 \, dx}{\left (
\int_{B_R} |x|^{-N}\, \left(-\log \left ( \frac{|x|}{R} \right)
\right)^{-\frac{2(N-1)}{N-2}}\, |v(|x|)|^{\frac{2N}{N-2}}\, dx
\right )^{\frac{N-2}{N}}},
\end{equation}
is
\begin{equation}\label{bestconstant}
C_{HS} := S(N)\, (N-2)^{-2(N-1)/N}
\end{equation}
and it is achieved by
\begin{equation} \label{extremalsw}
v_{\mu,\nu} (|x|) = \psi_{\mu,\nu} \left ( \left(-\log \left (
\frac{|x|}{R} \right) \right)^{-\frac{1}{N-2}} \right) = \left (
\mu^2 + \nu^2 \left(-\log \left ( \frac{|x|}{R} \right)
\right)^{-\frac{2}{N-2}} \right )^{-\frac{N-2}{2}},
\end{equation}
where the infimum is taken over the radially symmetric functions
of $W_{0}^{1,2} (|x|^{-(N-2)},B_R)$.
\end{theorem}
Observe that we may continuously define $u(0) = \mu^{-(N-2)}$ and
$u(R)=0$. We also have (see the proof of Theorem \ref{thmHSw}),
that $v_{\mu,\nu} \in W^{1,2}(|x|^{-(N-2)},B_R)$ but
$|x|^{-(N-2)}\, v_{\mu,\nu} \not\in H_0^1(B_R)$.

Concerning (\ref{HS1}), under the transformation (\ref{transx}),
we relate it with the Sobolev inequality (\ref{sobolev}). Then, we
prove that the best constant in (\ref{HS1}) is $C_{HS}$, as
defined in (\ref{bestconstant}) and the minimizers of (\ref{HS1}) are
\begin{equation}\label{minimizer}
u_{m,n} (|x|) = |x|^{-\frac{N-2}{2}} \psi_{m,n} \left( \left (-\log
\left ( \frac{|x|}{R} \right) \right)^{-\frac{1}{N-2}}
\right),\;\;\; x \in B_R \backslash \{0\},\;\;\; \phi_n
|_{\partial B_R} =0.
\end{equation}
where $\psi_{m,n}$ are given by (\ref{extremalsw});
\begin{theorem} \label{thm2}
The inequalities (\ref{HS1}) and
\begin{equation} \label{eq6}
\int_{\mathbb{R}^N} (\nabla w(t))^2\, dt \geq
(N-2)^{\frac{2(N-1)}{N}} C_{HS} \left ( \int_{\mathbb{R}^N}
|w(t)|^{\frac{2N}{N-2}}\, dt \right )^{\frac{N-2}{N}}.
\end{equation}
are equivalent under the transformation (\ref{transx}). Then,
the best constant in (\ref{HS1}) is (\ref{bestconstant}) and the
minimizers are given by (\ref{minimizer}).
\end{theorem}

In this direction, making some straightforward calculations, we have that
\begin{theorem} \label{thm1}
For each $n$, $\phi_n$ solves the corresponding to (\ref{HS1})
Euler-Lagrange equation:
\begin{eqnarray} \label{eqel1}
- \Delta u - \left ( \frac{N-2}{2} \right )^2 \frac{u}{|x|^2} &=&
\left(-\log \left ( \frac{|x|}{R} \right)
\right)^{-\frac{2(N-1)}{N-2}}\ u^{\frac{N+2}{N-2}} \\
u|_{\partial \Omega} &=& 0. \nonumber
\end{eqnarray}
\end{theorem}

For the nonradial case, i.e. $\Omega$ is an arbitrary bounded
domain in $\mathbb{R}^N$, containing the origin, $D_0 = \sup_{x
\in \Omega}$ and $D>D_0$, we refer to the recent work \cite{AFT}.
\vspace{0.1cm}

The paper is organized as follows: In Section 2 we consider the
spaces $W_{0}^{1,2} (|x|^{-(N-2)},\Omega)$ and $H(\Omega)$, we
prove Lemma \ref{lemmaWH} and as a consequence we obtain some
Caffarelli-Kohn-Nirenberg inequalities. In Section 3 we consider
inequality (\ref{m}) and in Section 3 we give the proof of the
remaining theorems.

For Hardy inequalities and their possible improvements we refer to
\cite{bm, bms, bv, cm2, davies, davdup03, delv04, ft02, ok, vz00}
and for various type of Hardy-Sobolev inequalities we refer to the
works \cite{ae05,
as,aft06,bt02,behl,cfms,ch03,fmt06,gk04,gr06,naz06,tt07,TZ07}.
\vspace{0.2cm} \\
{\bf Notation} In the sequel we often use the notation $r=|x|$.
%
%
\section{The spaces $W_{0}^{1,2}(|x|^{-(N-2)},\Omega)$ and
$H(\Omega)$}

\setcounter{equation}{0}
In this section we give some further properties for the spaces
$W_{0}^{1,2} (|x|^{-(N-2)},\Omega)$ and $H(\Omega)$ and give the
connection between them, i.e., we give the proof of Lemma
\ref{lemmaWH}.

Concerning $W_{0}^{1,2} (|x|^{-(N-2)},\Omega)$ from \cite[Lemma
2.1]{ft02} we have that
\begin{lemma} \label{ftW}
(i) If $u \in H_0^1 (\Omega)$, then $|x|^{\frac{N-2}{2}} u \in
W_{0}^{1,2} (|x|^{-(N-2)},\Omega)$. \\
(ii) If $v \in W_{0}^{1,2} (|x|^{-(N-2)},\Omega)$, then $|x|^{-a}
v \in H_0^1 (\Omega)$, for all $a<\frac{N-2}{2}$. \\
(iii) The norm
\[
\left( \int_{\Omega} |x|^{-(N-2)}\, |\nabla w|^2\, dx +
\int_{\Omega} |x|^{-(N-2)}\, w^2\, dx \right)^{1/2}
\]
is an equivalent norm for the space $W_{0}^{1,2}
(|x|^{-(N-2)},\Omega)$. \\
(iv) The space $W_{0}^{1,2} (|x|^{-(N-2)},\Omega)$ is a Hilbert
space with inner product
\[
<\phi,\psi>_{W_{0}^{1,2} (|x|^{-(N-2)},\Omega)} = \int_{\Omega}
|x|^{-(N-2)}\, \nabla \phi \cdot \nabla \psi\, dx.
\]
\end{lemma}

Concerning $H(\Omega)$ from \cite{vz00} we have that
\begin{lemma} \label{vzH}
(i) The space $H(\Omega)$ is a Hilbert space with inner product
\[
<\phi,\psi>_{H(\Omega)} = \int_{\Omega} \nabla \phi \cdot \nabla
\psi\, dx - \left ( \frac{N-2}{2} \right )^2\, \int_{\Omega}
\frac{\phi\, \psi}{|x|^2}\, dx - L(\phi)\, L(\psi).
\]
(ii) If $u \in H_0^1 (\Omega)$, then $u \in H(\Omega)$ and if $u
\in H(\Omega)$ then $u \in \cap_{q<2} W^{1,q}(\Omega)$, i.e.,
\[
H_0^1 (\Omega) \subset H(\Omega) \subset \cap_{q<2}
W^{1,q}(\Omega)
\]
(iii) The continuous imbedding $H(\Omega) \hookrightarrow
H_0^{s}$, $0 \leq s <1$ imply that the space $H(\Omega)$ is
compactly embedded in $L^{q} (\Omega)$, for any $1 \leq q <
\frac{2N}{N-2}$.
\end{lemma}
Moreover, from \cite[Theorem 4.2]{vz00} we have the following.
\begin{theorem} \label{Heigen}
Let $B_R$ the sphere in $\mathbb{R}^N$, $N \geq 3$, centered at
the origin with radius $R$. Let $z_{m,n}$ be the n-th zero of the
Bessel function $J_m$ and $\phi_k (\sigma)$ be the orthonormal
eigenfunctions of the Laplace-Beltrami operator with corresponding
eigenvalues $c_k = k(N+k-2)$, $k \geq 0$. Then, the two-parameter
family
\begin{equation}\label{vzeigen}
e_{k,n} (r,\sigma) = r^{-\frac{N-2}{2}}\, J_m \left (
\frac{z_{m,n}}{R} r \right )\, \phi_k(\sigma),
\end{equation}
with $m^2 = k(k+N-2)$, consist an orthogonal basis of $L^2 (B_R)$.
\end{theorem}
\begin{remark} \label{remH}
Note that all the $J_m$ vanish at $r=0$, except $J_0$ for which
(under normalization) $J_0(0) =1$. Then, the maximal singularity
corresponds to the sub-family of eigenfunctions with $j=0$
\[
e_{0,n} = O(r^{-\frac{N-2}{2}})
\]
These functions represent the complete sub-basis for the subspace
$X_1$ of radial functions in $L^2 (B_R)$. They do not belong to
$H_0^1 (\Omega)$ but belong to $H(\Omega)$.
\end{remark} \vspace{0.2cm}
{\bf Proof of Lemma \ref{lemmaWH}}\ \ (i) Let $v \in
C_{0}^{\infty} (\Omega)$, setting $u(x)= |x|^{-(N-2)/2}\, v(x)$ we
have
\begin{equation} \label{eqlwh1}
\int_{\Omega} |\nabla u|^2\, dx - \left ( \frac{N-2}{2} \right )^2
\int_{B_R} \frac{u^2}{|x|^2}\, dx = \int_{\Omega} |x|^{-(N-2)}
|\nabla v|^2\, dx + \frac{1}{2}\, \int_{\Omega} \nabla
|x|^{-(N-2)} \cdot \nabla  v^2\, dx.
\end{equation}
We first treat the radial case;  we assume that $\Omega = B_R$ and
let $v(r) \in C^{\infty} (0,R)$, $v(R) = 0$ and $v(r) \in
W_{0,r}^{1,2} (|x|^{-(N-2)},\Omega)$. From this point of view the
second integral in the right hand side of (\ref{eqlwh1}) is equal
to
\[
\frac{1}{2}\, \int_{\Omega} \nabla |x|^{-(N-2)} \cdot \nabla v^2\,
dx = - N\, \omega_N \frac{N-2}{2} \int_{0}^{R} (v^2)'\, dr =
\frac{N(N-2)}{2}\, \omega_N\, v^2(0).
\]
Then, (\ref{eqlwh1}) implies that $u \in H (B_R)$. For the
nonradial case, in order to estimate the second integral in the
right hand side of (\ref{eqlwh1}), we use the decomposition into
spherical harmonics; Let $v \in C_{0}^{\infty} (\Omega)$. If we
extend $u$ as zero outside $\Omega$, we may consider that $v \in
C_{0}^{\infty} (\mathbb{R}^N)$. Decomposing $v$ into spherical
harmonics we get
\[
v = \sum_{k=0}^{\infty} v_k := \sum_{k=0}^{\infty} f_k (r) \phi_k
(\sigma),
\]
where $\phi_k (\sigma)$ are the orthonormal eigenfunctions of the
Laplace-Beltrami operator with corresponding eigenvalues $c_k =
k(N+k-2)$, $k \geq 0$. The functions $f_k$ belong in
$C^{\infty}_{0} (\mathbb{R}^N)$, satisfying
\begin{equation}\label{fi}
f_k (r) = O (r^k),\;\;\; \mbox{and}\;\;\; f_k' (r) = O
(r^{k-1}),\;\;\; \mbox{as}\;\;\; r \downarrow 0.
\end{equation}
In particular, $\phi_0 (\sigma) =1$ and $v_0 (r) =
\frac{1}{|\partial B_r|} \int_{\partial B_r} u\, ds$, for any
$r>0$. Then, for any $k \in \mathbb{N}$, from (\ref{fi}) we have
that
\begin{eqnarray} \label{2.4}
\frac{1}{2}\, \int_{\Omega} \nabla |x|^{-(N-2)} \cdot \nabla v^2\,
dx &=& \frac{1}{2}\, \sum_{k=0}^{\infty} \int_{\mathbb{R}^N}
\nabla |x|^{-(N-2)} \cdot \nabla f_k^2\, dx \nonumber \\ &=&
\frac{N(N-2)}{2}\, \omega_N\, \sum_{k=0}^{\infty} \lim_{r \to 0}
f_k^2 (r) \nonumber \\ &=& \frac{N(N-2)}{2}\, \omega_N\,  f_0^2
(0) \nonumber \\ &=&  \frac{N(N-2)}{2}\, \omega_N\, v_0 (0).
\end{eqnarray}
Then, $u \in H(V)$ and (\ref{eqlwh1}) is equal to (\ref{normWH}).
\vspace{0.2cm}

(ii) Assume now that $u \in H(\Omega)$. Setting $v(x)=
|x|^{N-2/2}\, u(x)$ we have
\begin{eqnarray} \label{eqlwh2}
\int_{\Omega} |x|^{-(N-2)} |\nabla v|^2\, dx &=& \int_{\Omega}
|\nabla u|^2\, dx + \left ( \frac{N-2}{2} \right )^2 \int_{B_R}
\frac{u^2}{|x|^2}\, dx \nonumber \\ &&+ \frac{N-2}{2}\,
\int_{\Omega} |x|^{-2}\, x \cdot \nabla  u^2\, dx.
\end{eqnarray}
In order to estimate the last integral above, we use Theorem
\ref{Heigen}. Since $u \in H(\Omega)$ from Lemma \ref{vzH} we have
that $u \in L^2 (\Omega)$. Moreover, for some $R > \sup_{x \in
\Omega}$, we may assume that $u \in L^2 (B_R)$. Then, Theorem
\ref{Heigen} implies that there exist $c_n \in \mathbb{R}$,
$n=1,...$ such that
\[
u = \sum_{n=0}^{\infty} \sum_{k=0}^{\infty} c_n\,
\tilde{e}_{k,n}(r)\, \phi_k (\sigma),
\]
where
\[
\tilde{e}_{k,n}(r) := r^{-\frac{N-2}{2}}\, J_m \left (
\frac{z_{m,n}}{R} r \right ).
\]
Then,
\begin{eqnarray} \label{eqlwh3}
\frac{N-2}{2}\, \int_{\Omega} |x|^{-2}\, x \cdot \nabla u^2\, dx
&=& \frac{N-2}{2}\, \int_{B_R} |x|^{-2}\, x \cdot \nabla u^2\, dx
\nonumber \\ &=& \sum_{n=0}^{\infty} \sum_{k=0}^{\infty} c_n\,
\frac{N-2}{2}\, \int_{B_R} |x|^{-2}\, x \cdot \nabla
\tilde{e}_{k,n}^2(r)\, dx.
\end{eqnarray}
We have further, for every $k$ and $n$, that
\begin{eqnarray*}
\frac{N-2}{2}\, \int_{B_R} |x|^{-2}\, x \cdot \nabla
\tilde{e}_{k,n}^2(r)\, dx &=& \frac{N-2}{2}\, N\, \omega_N\,
\int_{0}^{R} r^{N-2}\, (\tilde{e}_{k,n}^2(r))'\, dr \\
&=& - \frac{(N-2)^2}{2}\, N\, \omega_N\, \int_{0}^{R} r^{N-3}\,
\tilde{e}_{k,n}^2(r)\, dr \\ && - \frac{N-2}{2}\, N\, \omega_N\,
\lim_{r \to 0} r^{N-2}\, \tilde{e}_{k,n}^2(r).
\end{eqnarray*}
or
\begin{eqnarray} \label{eqlwh4}
\sum_{n=0}^{\infty} \sum_{k=0}^{\infty} \frac{N-2}{2}\, c_n\,
\int_{B_R} \nabla |x|^{-2}\, x \cdot \nabla \tilde{e}_{k,n}^2(r)\,
dx = - \frac{(N-2)^2}{2}\, \int_{\Omega} \frac{u^2}{|x|^2}\, dx
\;\;\;\;\;\;\;\;\;\;\;\; \nonumber
\\ - \frac{N(N-2)}{2}\, \omega_N\, \sum_{n=0}^{\infty}
\sum_{k=0}^{\infty} c_n\,
 \lim_{r \to 0} r^{N-2}\, \tilde{e}_{k,n}^2(r).
\end{eqnarray}
However, Remark \ref{remH} implies that the only nonzero terms in
the above limit is given by $\tilde{e}_{0,n}^2(r)$. Hence,
(\ref{eqlwh3}) and (\ref{eqlwh4}) give that
\[
\frac{N-2}{2}\, \int_{\Omega} \nabla |x|^{-2}\, x \cdot \nabla
u^2\, dx = - \frac{(N-2)^2}{2}\, \int_{\Omega} \frac{u^2}{|x|^2}\,
dx  - \frac{N(N-2)}{2}\, \omega_N\, c_0.
\]
Finally, (\ref{eqlwh2}) becomes
\[
\int_{\Omega} |x|^{-(N-2)} |\nabla v|^2\, dx = \int_{\Omega}
|\nabla u|^2\, dx - \frac{(N-2)^2}{2}\, \int_{\Omega}
\frac{u^2}{|x|^2}\, dx  - \frac{N(N-2)}{2}\, \omega_N\, c_0.
\]
It is clear from the above discussion that $c_0$ corresponds to
$v_0 (0)$, i.e. we again derive (\ref{normWH}) and the proof is
completed.\ $\blacksquare$
\begin{corollary}
Assume now that $v_n$ is a bounded sequence in $W_{0}^{1,2}
(|x|^{-(N-2)},\Omega)$. Then $u_n = |x|^{-(N-2)/2}\, v_n$ is a
bounded sequence in $H(\Omega)$. The compact imbeddings of Lemma
\ref{vzH} imply that, up to some subsequence, $u_n$ converge in
$L^{q}(\Omega)$ to some $u$. Thus, we obtain the compact
imbeddings
\begin{equation}\label{compactimbed1}
W_{0}^{1,2} (|x|^{-(N-2)},\Omega) \hookrightarrow L^q
(|x|^{-q(N-2)/2},\Omega),\;\;\; \mbox{for any}\;\;\; 1 \leq q <
\frac{2N}{N-2}
\end{equation}
and since $1 \leq q$, we further obtain the compact imbeddings
\begin{equation}\label{compactimbed2}
W_{0}^{1,2} (|x|^{-(N-2)},\Omega) \hookrightarrow L^q
((|x|^{-(N-2)/2},\Omega),\;\;\; \mbox{for any}\;\;\; 1 \leq q <
\frac{2N}{N-2},
\end{equation}
where the weighted space $L^q (w(x),\Omega)$ is defined as the
closure of $C_{0}^{\infty}(\Omega)$ functions under the norm
\[
||\phi||_{L^q (w(x),\Omega)} = \left ( \int_{\Omega} w(x)\,
|\phi|^q\, dx \right )^{\frac{1}{q}}.
\]
\end{corollary}
\begin{remark}
In (\ref{compactimbed1}) is clear that $q$ cannot reach
$\frac{2N}{N-2}$. For this value of $q$ the best that we can have
is the inequality corresponding to (\ref{HSw}). In this sense the
results obtained in the previous Corollary complete the results
obtained in \cite{ckn} (see also \cite{cw01, cc93, ww03})
concerning the Caffarelli - Kohn - Nirenberg Inequalities, in the
limiting case where $a=\frac{N-2}{2}$.
\end{remark}
%
%
%
\section{Inequalities (\ref{m}) and (\ref{HS1})}
\setcounter{equation}{0}
In \cite[Theorem 4]{maz85} the following general inequality was
proved
\begin{equation} \label{m85}
\left [ \int_{-\infty}^{+\infty} \left | \int_{r}^{\infty} f(t)\,
dt \right |^q\, d\mu(r) \right ]^{1/q} \leq C\, \left [
\int_{-\infty}^{+\infty} |f(r)|^p\, d\nu(r)  \right ]^{1/p},
\end{equation}
where $1 \leq p \leq q \leq \infty$, which holds for any $f \in
C_{0}^{\infty}(\mathbb{R})$, if and only if the following quantity
\[
B = \sup_{l \in (-\infty,+\infty)} [\mu((-\infty,l))]^{1/q} \left
[ \int_{l}^{\infty} \left ( \frac{d\nu^*}{dr} \right
)^{-1/(p-1)}\, dr \right ]^{(p-1)/p},
\]
where $\nu^*$ is the absolutely continuous part of $\nu$, is
finite. Moreover, if $C$ is the best constant in (\ref{m85}), then
\begin{equation}\label{estbest}
B \leq C \leq B\, \left ( \frac{q}{q-1} \right )^{(p-1)/p}\,
q^{1/q}.
\end{equation}
Inequality (\ref{m}) is obtained (see \cite[Lemma
2.2]{ft02}) by setting $p=2$, $q=2N/(N-2)$, $f(r) = v'(r)$,
$\mu(r) = r^{-1}\, (-log r)^{-2(N-1)/(N-2)}\, \chi_{(0,1)} dr$
and $d\nu = r\, \chi_{(0,1)} dr$

As prof. V Maz'ya pointed to us inequality (\ref{m}) may also be
obtained from Bliss' inequality:
\begin{proposition} \label{pBliss}
For all $v:(0,\infty) \to \mathbb{R}$ absolutely continuous with
$v' \in L^k (0,\infty)$ and $v(0)=0$ one has
\begin{equation}\label{bliss}
\int_{0}^{\infty} \frac{|v|^l}{|x|^{l-h}}\, dx \leq K\, \left (
\int_{0}^{\infty} |v'|^k\, dx \right )^{l/k},
\end{equation}
where $l>k>1$, $h=l/k -1$ and
\[
K= \frac{1}{l-h-1}\, \left [ \frac{h\, \Gamma(l/h)}{\Gamma(1/h)\,
\Gamma((l-l)/h)} \right ]^h.
\]
Moreover, equality holds in (\ref{bliss}) if and only if
\begin{equation}\label{minbliss}
v(x) = (a+bx^{-h})^{-1/h},
\end{equation}
for arbitrary positive constants $a$ and $b$.
\end{proposition}
In the case now where $k=2$ and $l=2N/(N-2)$ we have that
$h=2/(N-2)$ and $l-h =2(N-1)/(N-2)$. Hence (\ref{bliss}) is equal
to
\begin{equation}\label{mbliss}
\int_{0}^{\infty} |v(t)|^{\frac{2N}{N-2}}\,
t^{-\frac{2(N-1)}{N-2}}\, dt \leq K\, \left ( \int_{0}^{\infty}
|v'(t)|^2\, dt \right )^{N/(N-2)}.
\end{equation}
As an alternative proof of Lemma \ref{lemmam} we may prove the
following.
\begin{lemma} \label{lemmamb}
(a)\ Inequality (\ref{mbliss}) under the change of variables
\[
t = -\log \left (\frac{r}{R} \right )
\]
is equivalent to (\ref{m})

(b)\ Inequality (\ref{mbliss}) under the change of variables
\[
t = r^{-(N-2)}
\]
is equivalent to (\ref{sobolev}).
\end{lemma}
We now give the proof of Lemma \ref{lemmam}. \vspace{0.2cm} \\
{\bf Proof of Lemma \ref{lemmam}}\ Let $v \in
C_{0}^{\infty}(0,R)$. Using the transformation (\ref{trans}) we
have that
\[
v'(r) = \frac{1}{N-2}\, w'(t)\, t^{N-1}\, r^{-1}
\]
and
\[
dr = (N-2)\, t^{-(N-1)}\, r.
\]
Then,
\[
\int_{0}^{R} r\, |v'|^2\, dr = \int_{0}^{\infty} r\,
\frac{1}{N-2}\, |w'(t)|^2\, t^{N-1}\, r^{-1}\, dt =\frac{1}{N-2}\,
\int_{0}^{\infty} t^{N-1}\, |w'(t)|^2\, dt
\]
and
\begin{eqnarray*}
\int_{0}^{R} r^{-1}\, \left ( - \log \left ( \frac{r}{R} \right )
\right )^{-\frac{2(N-1)}{N-2}} |v|^\frac{2N}{N-2}\, dr &=&
\int_{0}^{\infty} r^{-1}\, t^{2(N-1)}\, |w(t)|^\frac{2N}{N-2}\,
(N-2)\, t^{-(N-1)}\, r\, dt \\ &=& (N-2)\, \int_{0}^{\infty}
t^{N-1}\, |w(t)|^{\frac{2N}{N-2}}\, dt.
\end{eqnarray*}
So, inequality (\ref{m}) becomes
\[
\frac{1}{N-2}\, \int_{0}^{\infty} t^{N-1}\, |w'(t)|^2\, dt \geq c
(N-2)^{\frac{N-2}{N}} \left ( \int_{0}^{\infty} t^{N-1}\,
|w(t)|^{\frac{2N}{N-2}}\, dt  \right )^{\frac{N-2}{N}}
\]
or
\[
\int_{\mathbb{R}^N} |\nabla w|^2\, dx \geq c\,
(N-2)^{\frac{2(N-1)}{N}}\, (N\, \omega_N)^{\frac{2}{N}}\, \left (
\int_{\mathbb{R}^N} |w|^{\frac{2N}{N-2}}\, dx  \right
)^{\frac{N-2}{N}}.
\]
It is clear that if $v \in C_{0}^{\infty}(0,R)$ we have that $w
\in D^{1,2}(\mathbb{R}^N)$. Then, the best constant and the
minimizers are given by (\ref{bestm}) and (\ref{minimm}),
respectively and the proof is completed.\ $\blacksquare$
\vspace{0.2cm} \\
{\bf Proof of Lemma \ref{thm2}}\ Follows directly from Lemma
\ref{lemmam} and Lemma \ref{lemmaWH}. \ $\blacksquare$
\vspace{0.3cm} \\
{\bf Acknowledgements.} The author thanks prof. V. Maz'ya for
noticing the connection between inequality (\ref{m}) and the Bliss
inequality, prof. Luis Escauriaza for bringing to his attention
the work \cite{mus} and prof. A. Tertikas for informing, after a
personal communication, the results that were obtained in
\cite{AFT}.
\begin{small}
%
%
\bibliographystyle{amsplain}

\end{small}
\end{document}